\documentclass[letterpaper, 10 pt, conference]{ieeeconf}  

\IEEEoverridecommandlockouts                              

\usepackage{graphics}
\usepackage{graphicx}
\usepackage{epsfig}
\usepackage{mathptmx}
\usepackage{times}
\usepackage{amsmath}
\usepackage{amssymb}
\usepackage{cite}
\usepackage{algorithm}
\usepackage{caption}
\usepackage{subcaption}
\usepackage{algorithmic}
\usepackage{url}

\title{\LARGE \bf
A PDE-constrained Optimization Approach to Optimal Trajectory Planning under Uncertainty via Reflected Schrödinger Bridges
}

\author{Wenxin Liu$^{1}$ and Dante Kalise$^{2}$
\thanks{$^{1}$Wenxin Liu is with the Department of Mathematics,
        Imperial College London, London, UK
        {\tt\small wenxin.liu24@imperial.ac.uk}}%
\thanks{$^{2}$ Dante Kalise is with the Department of Mathematics,
        Imperial College London, London, UK
        {\tt\small dkaliseb@imperial.ac.uk}}%
}

\begin{document}

\maketitle
\thispagestyle{empty}
\pagestyle{empty}


\begin{abstract}
A computational PDE-constrained optimization approach is proposed for optimal trajectory planning under uncertainty by means of an associated Schr\"odinger Bridge Problem (SBP). The proposed SBP formulation is interpreted as the mean-field limit associated with the energy-optimal evolution of a particle governed by a stochastic differential equation (SDE) with nonlinear drift and reflecting boundary conditions, constrained to prescribed initial and terminal densities. The resulting mean-field system consists of a nonlinear Fokker-Planck equation coupled with a Hamilton-Jacobi-Bellman equation, subject to two-point boundary conditions in time and Neumann boundary conditions in space. Through the Hopf-Cole transformation, this nonlinear system is recast as a pair of forward-backward advection-diffusion equations, which are amenable to efficient numerical solution via a standard finite element discretization. The weak formulation naturally enforces reflecting boundary conditions without requiring explicit particle-boundary collision detection, thus circumventing the computational difficulties inherent to particle-based methods in complex geometries. Numerical experiments on challenging 3D maze configurations demonstrate fast convergence, mass conservation, and validate the optimal controls computed through reflected SDE simulations.
\end{abstract}
\section{Introduction}
Optimal pathfinding in constrained and complex environments under uncertainty remains a relevant challenge across robotics \cite{foderaro}, physics, and computational mathematics \cite{bkp}. Beyond algorithmic approaches, several physical-chemical systems have demonstrated that optimal paths can naturally emerge from gradient-driven transport phenomena. For instance, Lagzi \textit{et al.}~\cite{lagzi2010maze} and Lovass \textit{et al.}~\cite{lovass2015maze} showed that droplets or tracers can autonomously traverse mazes along the steepest descent of chemical or thermal potentials, effectively realizing physical analogues of optimal transport. These studies indicate that minimum energy trajectories can emerge from self-organization within constrained geometries under appropriate gradient fields.

Classical optimal transport (OT)~\cite{villani2021topics} provides a complete mathematical framework to describe the displacement of mass under prescribed costs and constraints. However, extending OT to domains with reflecting boundaries, drift fields, or dynamic obstacles remains analytically and computationally demanding. In robotics, for example, autonomous agents must plan motions in cluttered spaces while avoiding obstacles and following external or preferred drifts, conditions that naturally correspond to constrained OT formulations with reflection.

The Schr\"odinger Bridge problem (SBP) \cite{schrodinger1932theorie,chen2021stochastic} offers a stochastic perspective to optimal transport that naturally accounts for diffusion and uncertainty. In the dynamic formulation due to Benamou and Brenier \cite{benamou2000computational}, the SBP is posed as a PDE-constrained optimization problem: a kinetic energy functional is minimized over velocity fields subject to a Fokker-Planck equation evolving the density from prescribed initial to terminal distributions. Incorporating general prior drifts and reflecting boundary conditions, which are necessary for obstacle avoidance in confined domains, significantly complicates the problem. Existing approaches based on Wasserstein proximal recursion (WPR) \cite{caluya2019gradient,caluya2021wasserstein,caluya2021reflected} approximate PDE solutions via particle-based SDE simulations with interpolation. These methods become impractical on complex 3D geometries, and accurately enforcing reflection conditions for particle trajectories on curved, non-convex boundaries is notoriously difficult.

In this paper, we propose a PDE-based computational approach to the reflected SBP (RSBP) with nonlinear prior drift. Our approach departs from particle-based methods by working directly at the level of the PDE optimality system. The first-order conditions for the RSBP consist of a nonlinearly coupled forward Fokker-Planck and backward Hamilton-Jacobi-Bellman equation, both subject to Neumann boundary conditions that encode reflection. Via the Hopf-Cole transformation, this nonlinear system is exactly recast as a pair of linear forward-backward advection-diffusion equations, where the only nonlinearity is confined to the coupling through temporal boundary conditions and is resolved by a simple fixed-point iteration. The idea of linearizing the optimality system via Hopf-Cole has been explored in related mean-field game settings \cite{refId0}, but its use for the RSBP with prior drift and its combination with finite element (FEM) discretization on complex 3D geometries is new. The key advantages over particle-based methods are threefold: reflecting boundary conditions are enforced naturally through the weak formulation without requiring explicit particle-boundary collision detection; the FEM framework handles complex non-convex geometries in a principled way; and the fixed-point iteration converges in very few iterations, making the overall method computationally efficient.

The rest of the paper is organized as follows. Section~\ref{sec:rsbp} formulates the RSBP with prior drift, derives the optimality system, and introduces the Hopf-Cole transformation. Section~\ref{sec:fem} presents the FEM discretization and the fixed-point algorithm. Section~\ref{sec:numerics} demonstrates the method on 2D and 3D maze configurations, with a focus on numerical convergence and mass preservation features of the scheme.

\section{Reflected Schr\"odinger Bridge with Prior}\label{sec:rsbp}

\subsection{Problem Formulation}

Let $\Omega \subset \mathbb{R}^n$ be a smooth, bounded domain with boundary $\partial\Omega$. We consider the energy-optimal steering of a diffusive particle subject to reflecting boundary conditions and a divergence-free prior drift field that satisfies $\nabla\cdot v = 0 \text{ in } \Omega \text{ and } v\cdot n = 0 \text{ on } \partial\Omega$. The particle trajectory $X_t \in \Omega$ evolves according to the reflected SDE:
\begin{equation*}
dX_t = \big(v(X_t,t) + u(X_t,t)\big)\,dt + \sqrt{\varepsilon}\,dW_t + n(X_t)\,d\gamma_t,
\end{equation*}
where $v(x,t)$ is a given prior drift field, $u(x,t)$ is the control to be optimized, $n(x)$ is the inward unit normal at the boundary, and $\gamma_t$ is the boundary local time enforcing reflection. The particle satisfies $X_t \in \Omega$ for all $t \in [0,1]$, with prescribed initial and terminal distributions $X_0 \sim \rho_0$ and $X_1 \sim \rho_1$, and minimizes the expected kinetic energy cost
\begin{equation*}
\min_u \mathbb{E}\left[\int_0^1 \frac{1}{2}\|u(X_t,t)\|^2 \,dt\right]
\end{equation*}
while ensuring the particle's distribution evolves from $\rho_0$ to $\rho_1$. We now make the particle problem precise at the population level.

\subsection{Mean-Field Formulation}

The probability density $\rho(x,t)$ of the particle satisfies a Fokker-Planck equation with reflecting boundary conditions. The reflected Schr\"odinger Bridge Problem (RSBP) corresponds to minimizing the kinetic energy
\begin{equation*}
J(\rho,u):=\frac{1}{2}\int_0^1\!\!\int_\Omega \rho(x,t)\|u(x,t)\|^2 \,dx\,dt
\end{equation*}
subject to the Fokker-Planck equation
\begin{equation*}
\partial_t \rho = \frac{\varepsilon}{2}\Delta\rho - \nabla\cdot\big(\rho(v + u)\big) \quad \text{in } \Omega \times (0,1),
\end{equation*}
with Neumann boundary conditions enforcing reflection
\begin{equation}\label{eq:neumann-rho}
\quad n \cdot \left(\frac{\varepsilon}{2}\nabla\rho - \rho(\nabla\lambda + v)\right) = 0  \quad \text{on } \partial\Omega \times (0,1),
\end{equation}
and temporal boundary conditions
\begin{equation}\label{eq:temporal-bc}
\rho(x,0) = \rho_0(x), \quad \rho(x,1) = \rho_1(x).
\end{equation}

This optimization problem seeks the optimal velocity field $u^*(x,t)$ that transports the density from $\rho_0$ to $\rho_1$ while minimizing kinetic energy and respecting geometric confinement. We derive its optimality conditions next.

\subsection{First-order Optimality Conditions}

The optimality system for the reflected Schr\"odinger bridge problem follows from a standard Lagrangian approach, see e.g. \cite{mfoc,bkp}, and consists of a Fokker-Planck equation coupled with a Hamilton-Jacobi-Bellman equation. Introducing the adjoint variable $\lambda(x,t)$, the optimality conditions for $\rho$ and $\lambda$ are given by:
\begin{align}
\partial_t \rho + \nabla\cdot\big(\rho(\nabla\lambda + v)\big) &= \frac{\varepsilon}{2}\Delta\rho, \label{eq:FP-optimality}\\
\partial_t \lambda + \frac{1}{2}\|\nabla\lambda\|^2 + \nabla\lambda \cdot v &= -\frac{\varepsilon}{2}\Delta\lambda, \label{eq:HJB-optimality}
\end{align}
subject to reflecting \eqref{eq:neumann-rho} and Neumann boundary conditions for the adjoint
\begin{equation*}
\nabla\lambda \cdot n = 0, \quad \text{on } \partial\Omega \times (0,1),
\end{equation*}
and temporal boundary conditions \eqref{eq:temporal-bc}.

The Fokker-Planck equation~\eqref{eq:FP-optimality} evolves forward in time from the initial condition, while the Hamilton-Jacobi-Bellman equation~\eqref{eq:HJB-optimality} is solved backward from a terminal condition determined by the requirement that $\rho(x,1) = \rho_1(x)$. This two-point boundary value problem in time, coupled with the nonlinear structure of both PDEs, presents significant computational challenges. Once this system has been solved, the optimal control is given by $u^*(x,t) = \nabla\lambda(x,t)$. The main computational difficulty is the nonlinear coupling between~\eqref{eq:FP-optimality} and~\eqref{eq:HJB-optimality}; we resolve it next via a change of variables.

\subsection{Hopf-Cole Transformation}

The Hopf-Cole transformation provides a fundamental linearization of the coupled nonlinear system~\eqref{eq:FP-optimality}--\eqref{eq:HJB-optimality}. We introduce potential functions $\varphi(x,t)$ and $\hat{\varphi}(x,t)$ through the change of variables:
\begin{equation}\label{eq:hopf-cole-def}
\varphi(x,t) = e^{\lambda(x,t)/\varepsilon}, \quad \hat{\varphi}(x,t) = \rho(x,t)e^{-\lambda(x,t)/\varepsilon}.
\end{equation}
The density factorizes as $\rho(x,t) = \varphi(x,t)\hat{\varphi}(x,t)$, and the optimal control can be recovered from:
\begin{equation}\label{eq:control-recovery}
u^*(x,t) = \nabla\lambda(x,t) = \varepsilon\left(\frac{\nabla\varphi}{\varphi}\right).
\end{equation}
Assuming that the prior drift satisfies
\[
\nabla\cdot v = 0 \quad \text{in }\Omega,
\qquad
v\cdot n = 0 \quad \text{on }\partial\Omega,
\]
we substitute the transformation~\eqref{eq:hopf-cole-def} into the optimality system~\eqref{eq:FP-optimality}--\eqref{eq:HJB-optimality}, from which we obtain that the potentials $\varphi$ and $\hat{\varphi}$ satisfy the forward-backward system:
\begin{align}
\partial_t \varphi = -\nabla\varphi \cdot v - \frac{\varepsilon}{2}\Delta\varphi,\qquad
\partial_t \hat{\varphi} = -\nabla\hat{\varphi} \cdot v + \frac{\varepsilon}{2}\Delta\hat{\varphi}\,, \label{eq:backward-phi}
\end{align}
together with homogeneous Neumann boundary conditions
\begin{equation}\label{eq:neumann-potentials}
\nabla\varphi \cdot n = 0,
\qquad
\nabla\hat{\varphi} \cdot n = 0
\quad \text{on }\partial\Omega\times (0,1)\,.
\end{equation}
The system is coupled by the endpoint conditions
\begin{equation}\label{eq:temporal-potentials}
\varphi(x,0)\,\hat{\varphi}(x,0) = \rho_0(x),
\qquad
\varphi(x,1)\,\hat{\varphi}(x,1) = \rho_1(x).
\end{equation}

The remarkable feature of this transformation is that equations~\eqref{eq:backward-phi} are linear advection-diffusion equations. The nonlinearity of the original system is now entirely contained in the coupling through the temporal boundary conditions~\eqref{eq:temporal-potentials}, which can be resolved via fixed-point iteration. Note that when there is no prior drift ($v \equiv 0$), the system~\eqref{eq:backward-phi} reduces to pure heat equations:
\begin{align*}
\partial_t \varphi = -\frac{\varepsilon}{2}\Delta\varphi, \qquad
\partial_t \hat{\varphi} = \frac{\varepsilon}{2}\Delta\hat{\varphi},
\end{align*}
which are linear and amenable to efficient numerical solution via standard parabolic PDE discretization methods. Section~\ref{sec:fem} describes the chosen discretisation.

\section{Finite Element Discretization}\label{sec:fem}

\subsection{Spatial Discretization}

We discretize the spatial domain $\Omega \subset \mathbb{R}^d$ using a conforming triangulation $\mathcal{T}_h = \{K\}$ with characteristic mesh size $h$. The finite element space consists of continuous piecewise linear functions:
\begin{equation*}
V_h := \{v_h \in C^0(\Omega) : v_h|_K \in \mathbb{P}_1(K) \;\forall K \in \mathcal{T}_h\},
\end{equation*}
where $\mathbb{P}_1(K)$ denotes the space of linear polynomials on element $K$. Let $\{\phi_i\}_{i=1}^{N_h}$ denote the nodal basis functions satisfying $\phi_i(x_j) = \delta_{ij}$, where $\{x_i\}_{i=1}^{N_h}$ are the mesh vertices. Any function $w_h \in V_h$ admits the representation
\begin{equation*}
w_h(x) = \sum_{i=1}^{N_h} w_i \phi_i(x),
\end{equation*}
where $w_i = w_h(x_i)$ correspond to the nodal values. Time is discretised as follows.

\subsection{Temporal Discretization}

We partition the time interval $[0,1]$ uniformly with time step $\Delta t = 1/K$, defining discrete times $t_k = k\Delta t$ for $k = 0,1,\ldots,K$. For temporal discretization, we employ the backward Euler method, which provides unconditional stability and is crucial for maintaining positivity of the density throughout the evolution. The space-time discretization of the forward-backward system~\eqref{eq:backward-phi} is given by the following weak formulations.

\subsection{Weak Formulation}
The weak formulation of the forward equation in ~\eqref{eq:backward-phi} at time level $k+1$ reads: find $\hat{\varphi}_h^{k+1} \in V_h$ such that for all test functions $w_h \in V_h$:
\begin{align}
\int_\Omega \frac{\hat{\varphi}_h^{k+1} - \hat{\varphi}_h^k}{\Delta t} w_h  + (\nabla\hat{\varphi}_h^{k+1} \cdot v) w_h + \frac{\varepsilon}{2}\nabla\hat{\varphi}_h^{k+1} \cdot \nabla w_h \,\,dx = 0.\label{eq:weak-forward}
\end{align}
The homogeneous Neumann boundary condition~\eqref{eq:neumann-potentials} is naturally incorporated through the weak formulation. After integration by parts of the diffusion term, the boundary integral vanishes due to $\nabla\hat{\varphi} \cdot n = 0$, eliminating the need for explicit boundary condition enforcement.

Similarly, the weak formulation of the backward equation at time level $k$, solved backward from $t = 1$ to $t = 0$, reads: find $\varphi_h^k \in V_h$ such that for all $w_h \in V_h$:
\begin{equation}
\int_\Omega \frac{\varphi_h^k - \varphi_h^{k+1}}{\Delta t} w_h+  (\nabla\varphi_h^k \cdot v) w_h  + \frac{\varepsilon}{2} \nabla\varphi_h^k \cdot \nabla w_h \,dx = 0.\label{eq:weak-backward}
\end{equation}

Again, the reflecting boundary conditions are naturally satisfied in the weak formulation. Substituting the nodal expansion yields the discrete systems to be solved at each time step.

\subsection{Discrete Linear Systems}
Expanding the discrete solutions in the nodal basis:
\begin{equation*}
\hat{\varphi}_h^k(x) = \sum_{i=1}^{N_h} \hat{\Phi}_i^k \phi_i(x), \quad \varphi_h^k(x) = \sum_{i=1}^{N_h} \Phi_i^k \phi_i(x),
\end{equation*}
and substituting into the weak formulations~\eqref{eq:weak-forward}--\eqref{eq:weak-backward}, we obtain linear systems of the form:
\begin{align}
\left(M + \Delta t C[v^k] + \frac{\varepsilon}{2}\Delta t L\right)\hat{\boldsymbol{\Phi}}^{k+1} &= M\hat{\boldsymbol{\Phi}}^k, \label{eq:forward-system}\\
\left(M + \Delta t C[v^k] + \frac{\varepsilon}{2}\Delta t L\right)\boldsymbol{\Phi}^k &= M\boldsymbol{\Phi}^{k+1}, \label{eq:backward-system}
\end{align}
where $\hat{\boldsymbol{\Phi}}^k = (\hat{\Phi}_1^k,\ldots,\hat{\Phi}_{N_h}^k)^\top$ and $\boldsymbol{\Phi}^k = (\Phi_1^k,\ldots,\Phi_{N_h}^k)^\top$ are the vectors of nodal values. The matrices are defined as:
\begin{align*}
M_{ij}&= \int_\Omega \phi_i \phi_j \,dx, \quad
L_{ij}= \int_\Omega \nabla\phi_i \cdot \nabla\phi_j \,dx,\\
C_{ij}[v]& = \int_\Omega (\nabla\phi_j \cdot v) \phi_i \,dx.
\end{align*}

All three matrices are sparse, symmetric (for $M$ and $L$), and can be efficiently assembled using standard finite element routines. For the pure diffusion case ($v \equiv 0$), the convection matrix vanishes ($C \equiv 0$), further simplifying the system. At each time step, we solve the linear systems~\eqref{eq:forward-system}--\eqref{eq:backward-system} using an iterative method. The system matrix:
\begin{equation*}
A = M + \Delta t C[v] + \frac{\varepsilon}{2}\Delta t L
\end{equation*}
is sparse and, for sufficiently small convection, symmetric positive definite. We employ the Generalized Minimal Residual (GMRES) method with incomplete LU (ILU) preconditioning to solve these systems efficiently. The preconditioner significantly accelerates convergence, particularly for refined meshes where the system becomes increasingly ill-conditioned \cite{precmfg}. It remains to address the coupling between forward and backward solves.

\subsection{Fixed-Point Algorithm}
We now address the solution of the coupled system via a fixed-point iteration. A critical implementation detail is ensuring that the potentials $\varphi_h$ and $\hat{\varphi}_h$ remain strictly positive throughout the computation. This is essential for computing the optimal control via~\eqref{eq:control-recovery}, which involves division by $\varphi$ and $\hat{\varphi}$; maintaining physical consistency of the density $\rho_h = \varphi_h \hat{\varphi}_h$; and ensuring numerical stability of the fixed-point iteration. After solving each linear system, we apply the projection
\begin{equation*}
\varphi_h \leftarrow \max(\varphi_h, \delta), \quad \hat{\varphi}_h \leftarrow \max(\hat{\varphi}_h, \delta),
\end{equation*}
where $\delta \approx 10^{-12}$ is a small positive threshold to account for round-off errors.

The coupling between the forward and backward equations through the temporal boundary conditions~\eqref{eq:temporal-potentials} creates a two-point boundary value problem. We resolve this coupling via a fixed-point iteration detailed in Algorithm \ref{alg:fixed-point}.

\begin{algorithm}[!ht]
\caption{Fixed-Point Iteration for RSBP}
\label{alg:fixed-point}
\begin{algorithmic}[1]
\STATE Initialize $\hat{\varphi}_h^{(0)}(x,T) \equiv 1$, $m \gets 0$, and $\text{error} \gets \infty$
\WHILE{$\text{error} > \text{tol}$}
    \STATE Compute terminal condition: \[\varphi_h^{(m)}(x,T) = \rho_1(x)/\hat{\varphi}_h^{(m)}(x,T)\]
    \STATE Solve backward equation~\eqref{eq:weak-backward} for $\varphi_h^{(m)}$
    \STATE Apply positivity projection $\varphi_h^{(m)}=\max(\varphi_h^{(m)}, \delta)$
    \STATE Compute initial condition: \[\hat{\varphi}_h^{(m+1)}(x,0) = \rho_0(x)/\varphi_h^{(m)}(x,0)\]
    \STATE Solve forward equation~\eqref{eq:weak-forward} for $\hat{\varphi}_h^{(m+1)}$
    \STATE Apply positivity projection $\hat{\varphi}_h^{(m+1)}=\max(\hat{\varphi}_h^{(m+1)}, \delta)$
    \STATE Compute density: $\rho_h^{(m+1)} = \varphi_h^{(m)} \hat{\varphi}_h^{(m+1)}$
    \STATE $\text{error} \gets \|\hat{\varphi}_h^{(m+1)}(T) - \hat{\varphi}_h^{(m)}(T)\|_{L^2}$
    \STATE $m \gets m + 1$
\ENDWHILE
\end{algorithmic}
\end{algorithm}

The convergence criterion monitors the change in $\hat{\varphi}$ at the terminal time, which directly controls satisfaction of the terminal density constraint. The fixed-point map $\hat{\varphi}^{(m)}(1)$ to $\hat{\varphi}^{(m+1)}(1)$ is contractive in Hilbert's projective metric, ensuring geometric convergence to the unique solution when the potentials remain strictly positive and bounded \cite{CGP_hilbert}. Standard finite element theory provides error estimates for the spatial and temporal discretization \cite{Quarteroni1994}. For P1 elements and backward Euler time-stepping, we have:
\begin{equation*}
\|\rho(t_n) - \rho_h^n\|_{L^2(\Omega)} \leq C(h^2 + \Delta t),
\end{equation*}
where the constant $C$ depends on model parameters and is independent of $h$ and $\Delta t$. We now demonstrate the method on non-trivial two- and three-dimensional geometries.
\section{Numerical experiments}\label{sec:numerics}
All experiments were implemented in FEniCS with a PETSc backend, using GMRES preconditioned with BoomerAMG (\texttt{hypre\_amg}) for the linear solves. Computations were run on a Linux workstation with dual Intel Xeon E5-2673 v4 processors (40 physical cores) and 251\,GiB RAM, using 16 MPI processes.

\subsection{Convergence study on a 2D maze}

We begin with a 2D test that permits a systematic convergence analysis. The domain $\Omega \subset [0,1]^2$ is a maze geometry consisting of a unit square with cross-shaped interior obstacles. Initial and terminal densities are Gaussians placed in opposite corners of the maze ($\sigma = 0.05$). The diffusion coefficient is $\varepsilon = 0.05$.

A divergence-free field $v$ with $\nabla\cdot v = 0$ in $\Omega$ and $v\cdot n = 0$ on $\partial\Omega$ is obtained by projecting a physically motivated velocity $V$ onto the divergence-free subspace via a Helmholtz--Hodge decomposition. Specifically, we solve the Poisson equation
\[
-\Delta \psi = \nabla\cdot V \quad\text{in }\Omega, \qquad
\frac{\partial\psi}{\partial n} = 0 \quad\text{on }\partial\Omega,
\]
and set $v = V - \nabla\psi$. This removes the irrotational component of $V$ and yields a numerically divergence-free velocity tangent to the reflecting boundary.

\begin{figure}[!ht]
\centering
\includegraphics[width=\columnwidth]{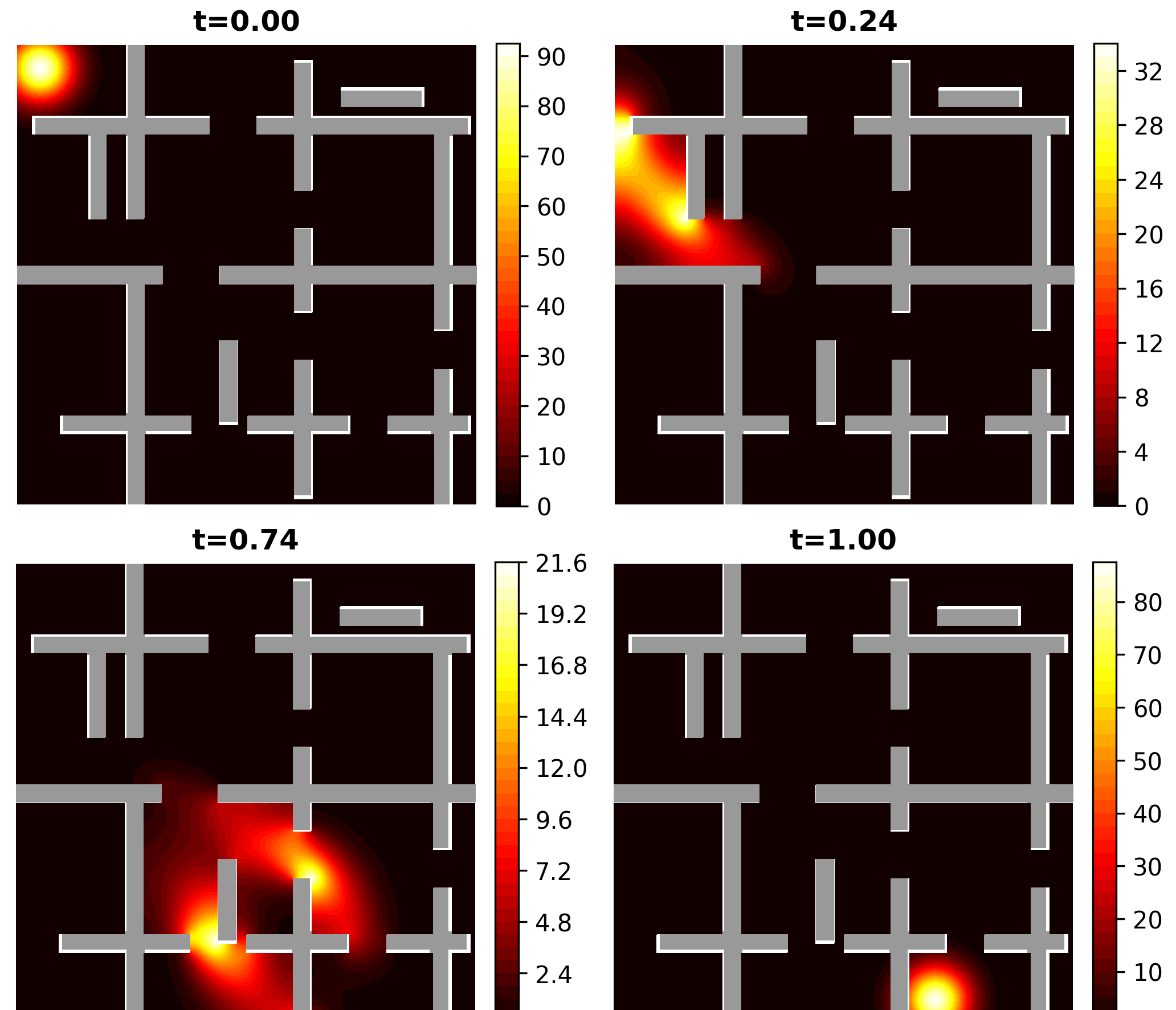}
\caption{2D maze, RSBP with divergence-free prior drift: the density $\rho(t,x)$ navigates the interior obstacles and arrives at the target corner at $t=1$.}
\label{fig:density_evolution_2d}
\end{figure}

Figure~\ref{fig:density_evolution_2d} shows the density evolving from the upper-left entrance to the lower-right exit under the divergence-free prior drift, confirming that the computed transport navigates the maze walls correctly.

\paragraph{Spatial and temporal convergence.}
We assess convergence of the transport cost $J$ and mass conservation error $\|\delta m\|_\infty = \max_t |\int_\Omega \rho\,dx - 1|$ under systematic mesh and time-step refinement, using the finest available resolution as reference ($N=2000$, $K=1200$). Tables~\ref{tab:spatial_2d} and~\ref{tab:temporal_2d} report the energy and mass convergence under spatial (fixed $K=1200$) and temporal (fixed $N=2000$) refinement, respectively.

Spatial convergence of the energy reaches order $\mathcal{O}(h^{1.87})$ in the pre-saturation range $N \in [100,500]$; for $N \geq 1000$ the error saturates at the level of the fixed-point stopping tolerance. The mass error converges at $\mathcal{O}(h^{2.15})$ without saturation, reaching $\mathcal{O}(10^{-5})$ at $N=2000$. Temporal convergence of the energy follows $\mathcal{O}(\Delta t^{1.58})$. Mass conservation under temporal refinement is insensitive to $K$ at the $10^{-5}$ level, confirming that mass errors are dominated by spatial resolution.

\begin{table}[!ht]
\centering
\caption{2D maze: spatial convergence of energy and mass error ($K=1200$).}
\label{tab:spatial_2d}
\renewcommand{\arraystretch}{1.05}
\begin{tabular}{rrccc}
\hline
$N$ & DOFs & $J_h$ & $|J_h - J_{\text{ref}}|$ & $\|\delta m\|_\infty$ \\
\hline
100  & 8{,}193       & 1.041544 & $3.08\times10^{-2}$ & $7.64\times10^{-3}$ \\
200  & 32{,}345      & 1.019818 & $9.04\times10^{-3}$ & $1.75\times10^{-3}$ \\
500  & 200{,}298     & 1.012313 & $1.54\times10^{-3}$ & $2.60\times10^{-4}$ \\
1000 & 795{,}748     & 1.012211 & $1.44\times10^{-3}$ & $6.16\times10^{-5}$ \\
2000 & 3{,}177{,}965 & 1.010775 & -- (ref.)           & $1.14\times10^{-5}$ \\
\hline
\end{tabular}
\end{table}

\begin{table}[!ht]
\centering
\caption{2D maze: temporal convergence of energy and mass error ($N=2000$).}
\label{tab:temporal_2d}
\renewcommand{\arraystretch}{1.05}
\begin{tabular}{rccc}
\hline
$K$ & $\Delta t$ & $|J_K - J_{\text{ref}}|$ & $\|\delta m\|_\infty$ \\
\hline
100  & $10^{-2}$        & $1.98\times10^{-1}$ & $9.73\times10^{-6}$ \\
200  & $5\times10^{-3}$ & $1.03\times10^{-1}$ & $1.39\times10^{-5}$ \\
500  & $2\times10^{-3}$ & $3.17\times10^{-2}$ & $1.28\times10^{-5}$ \\
1000 & $10^{-3}$        & $4.69\times10^{-3}$ & $1.50\times10^{-5}$ \\
1200 & $8\times10^{-4}$ & -- (ref.)           & $1.14\times10^{-5}$ \\
\hline
\end{tabular}
\end{table}

\subsection{Path planning in a 3D spiral maze}
\begin{figure}[!ht]
    \centering
    \begin{subfigure}{0.7\columnwidth}
        \centering
        \includegraphics[width=\linewidth]{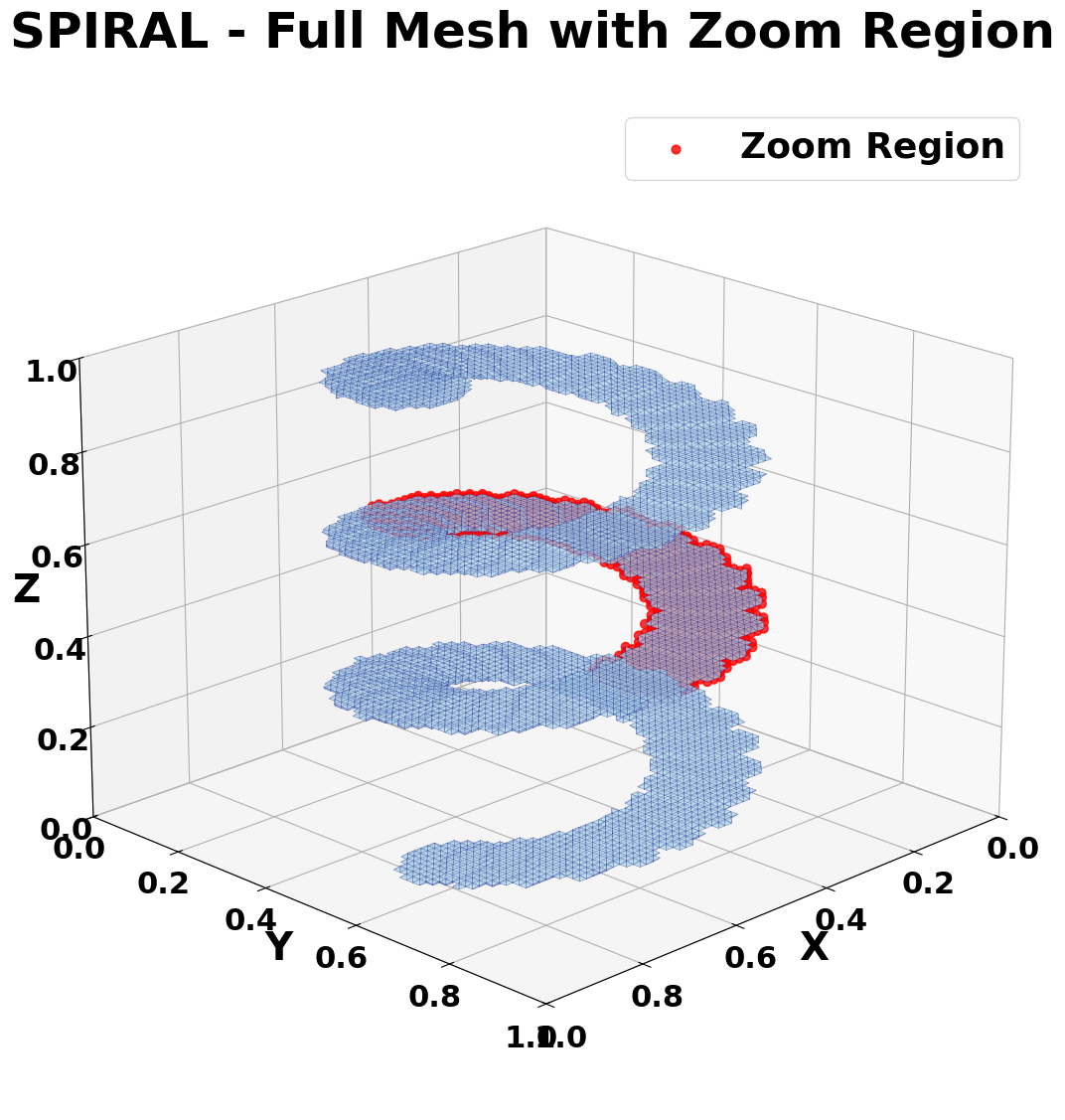}
        \caption{Full mesh}
    \end{subfigure}\\
    \begin{subfigure}{0.7\columnwidth}
        \centering
        \includegraphics[width=\linewidth]{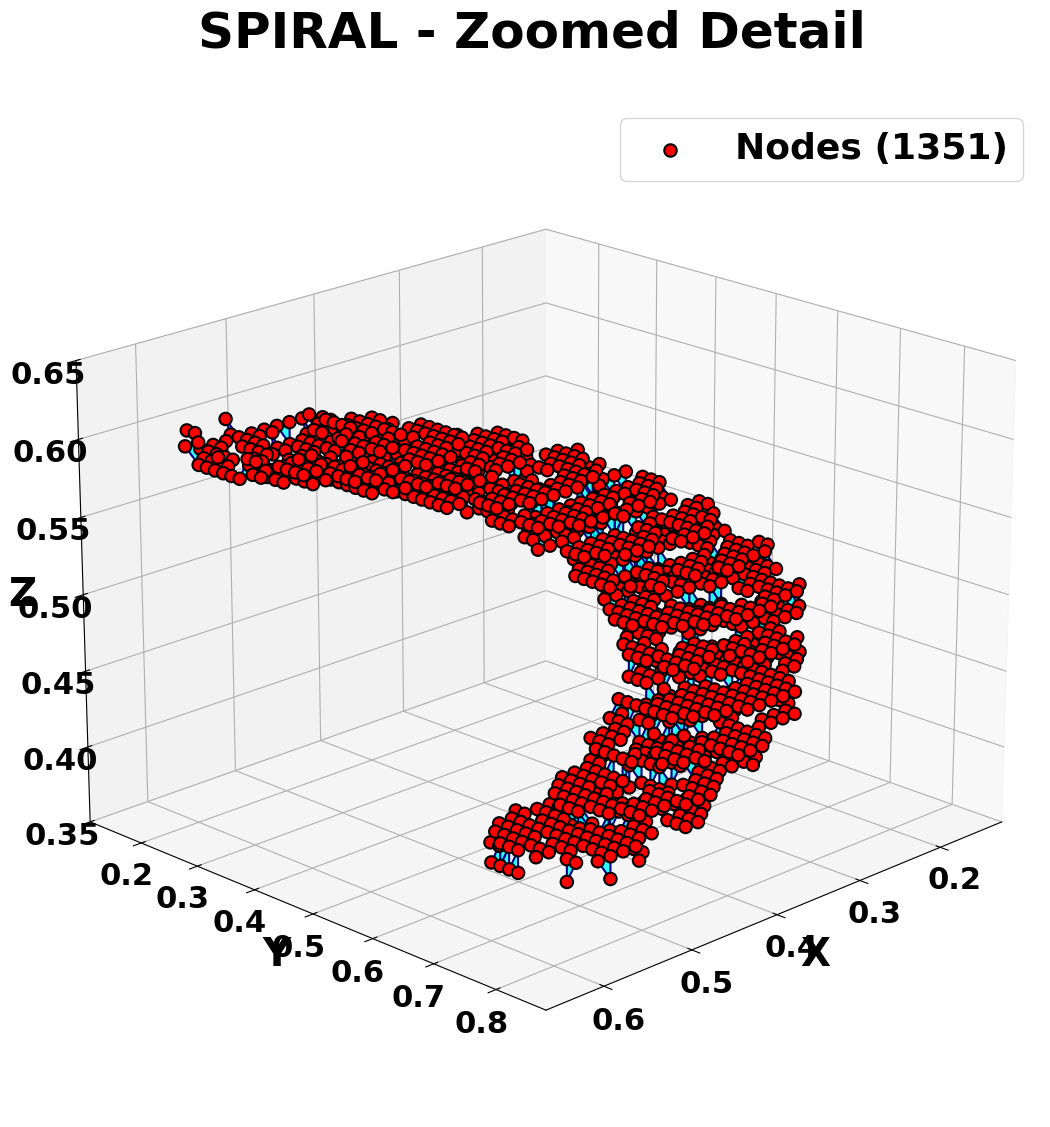}
        \caption{Zoomed mesh detail}
    \end{subfigure}
\caption{3D helical mesh and zoomed detail.}\label{fig:helical_mesh}
\end{figure}

We assess the proposed approach on a three-dimensional domain consisting of a helical tube parametrised as
\begin{equation*}
c(z) = \left( 0.5 + 0.25 \cos(6\pi z), \, 0.5 + 0.25 \sin(6\pi z), \, z \right),
\end{equation*}
for $z\in [0,1]$ with tube radius $r_t=0.1$. Figure~\ref{fig:helical_mesh} shows the finite element mesh of conforming P1 elements, with spatial resolutions ranging from $N = 60$ to $N = 90$ nodes per axis. Initial and terminal densities are given by
\begin{align*}
\rho_i(x) &= \frac{1}{Z_0} \exp\left(-\frac{\|x - x_i\|^2}{2\sigma^2}\right),\quad i=\{0,1\},
\end{align*}
where $x_0 = (0.72, 0.63, 0.05)$ (entrance, bottom), $x_1 = (0.76, 0.42, 0.95)$ (exit, top), and $\sigma = 0.05$. We set diffusion coefficient $\varepsilon = 0.5$, and fixed-point convergence relative tolerance to $10^{-2}$.

\begin{figure}[!ht]
    \centering
    \begin{subfigure}{0.9\columnwidth}
        \centering
        \includegraphics[width=\linewidth]{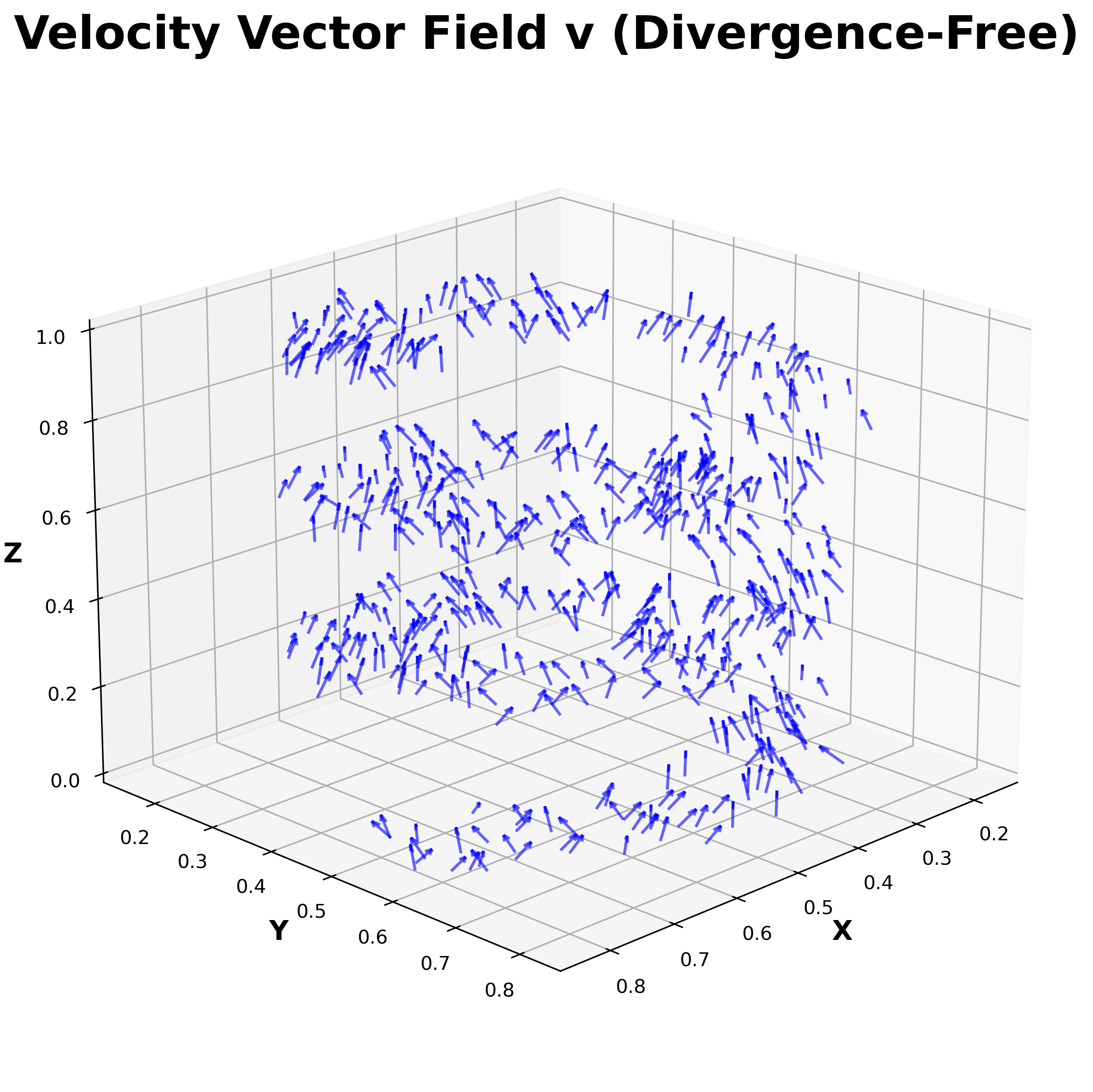}
        \caption{Prior drift field $v$}
    \end{subfigure}\\[4pt]
    \begin{subfigure}{0.9\columnwidth}
        \centering
        \includegraphics[width=\linewidth]{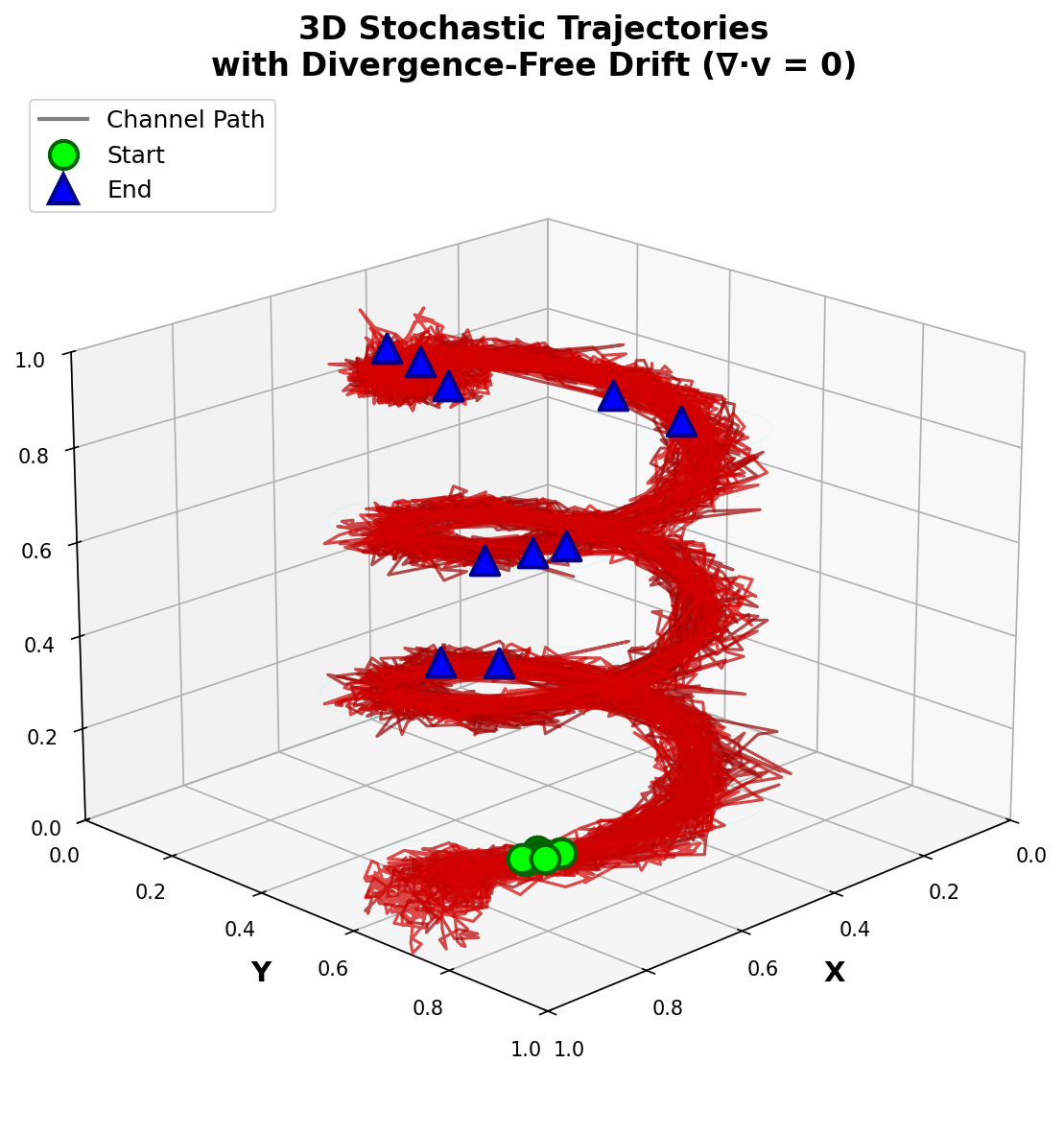}
        \caption{3D controlled trajectories}
    \end{subfigure}
    \caption{RSBP with prior drift: prior drift field and reflected SDE validation with 10 particles evolving under the optimal feedback~\eqref{sde_part}.}
    \label{fig:drift_particle_trajectories}
\end{figure}

\begin{figure}[!ht]
    \centering
    \begin{subfigure}{0.58\columnwidth}
        \centering
        \includegraphics[width=\linewidth]{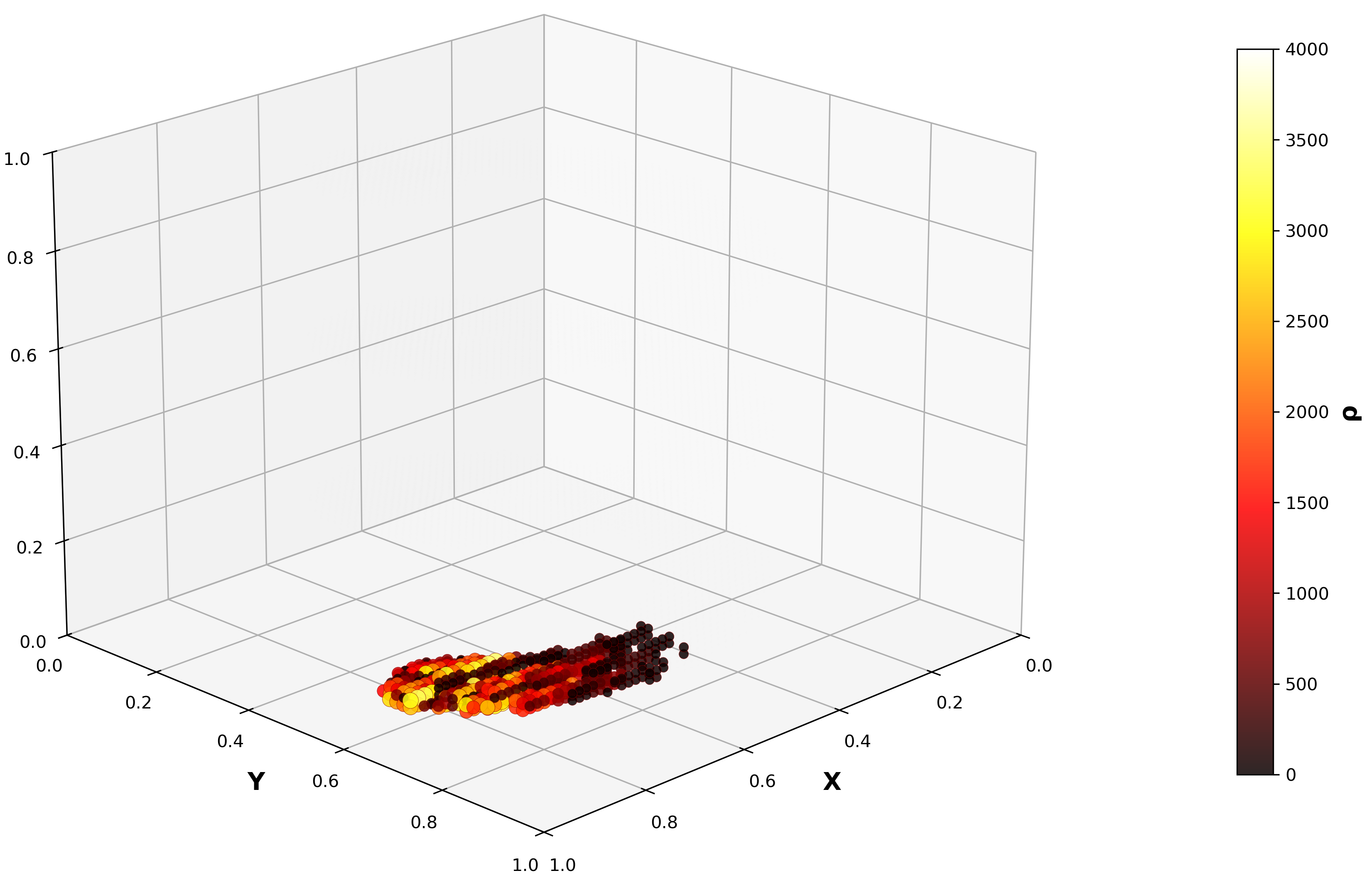}
        \caption{$t=0.0$}
    \end{subfigure}\hfill
    \begin{subfigure}{0.39\columnwidth}
        \centering
        \includegraphics[width=\linewidth]{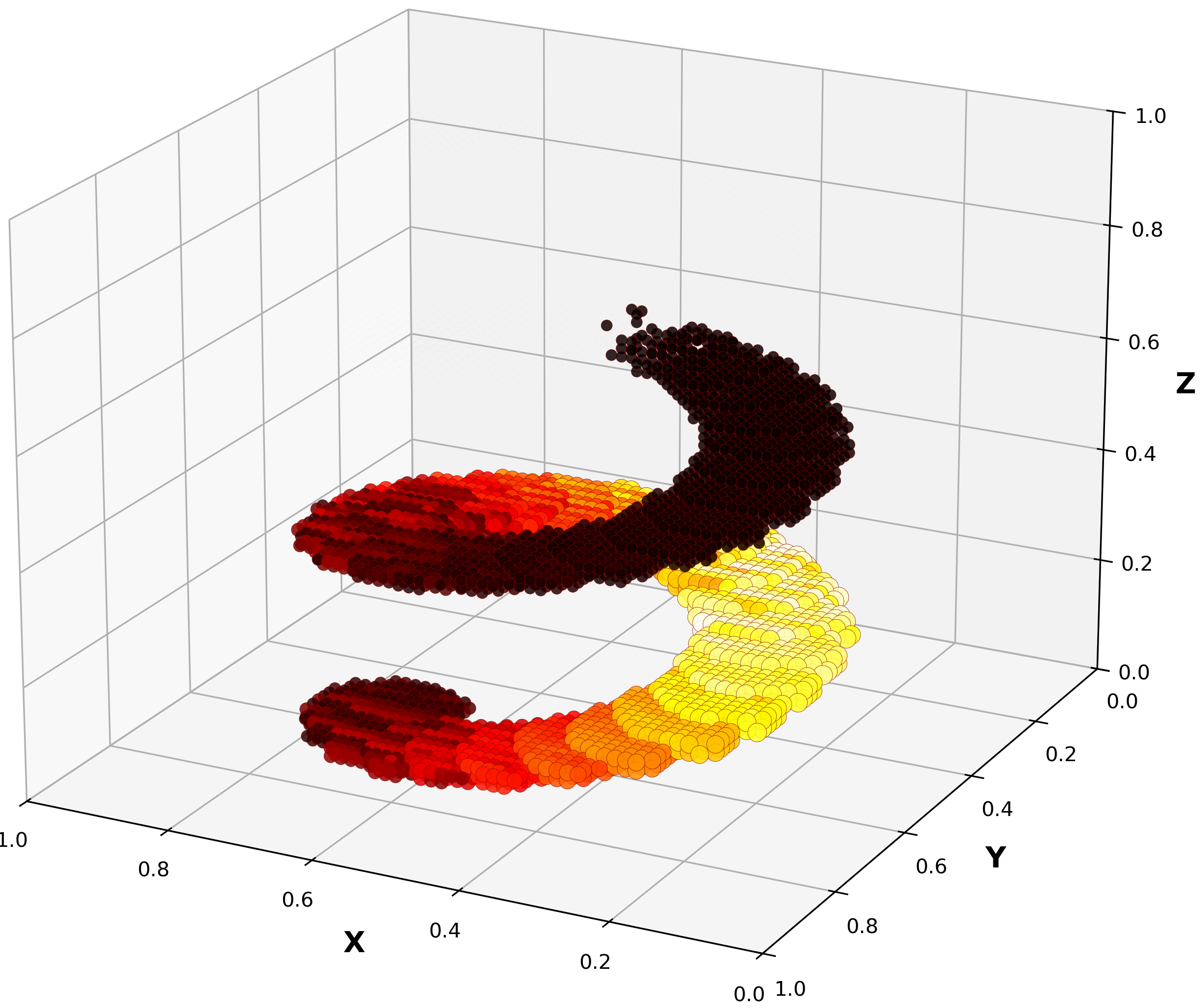}
        \caption{$t=0.2$}
    \end{subfigure}\\[4pt]
    \begin{subfigure}{0.48\columnwidth}
        \centering
        \includegraphics[width=\linewidth]{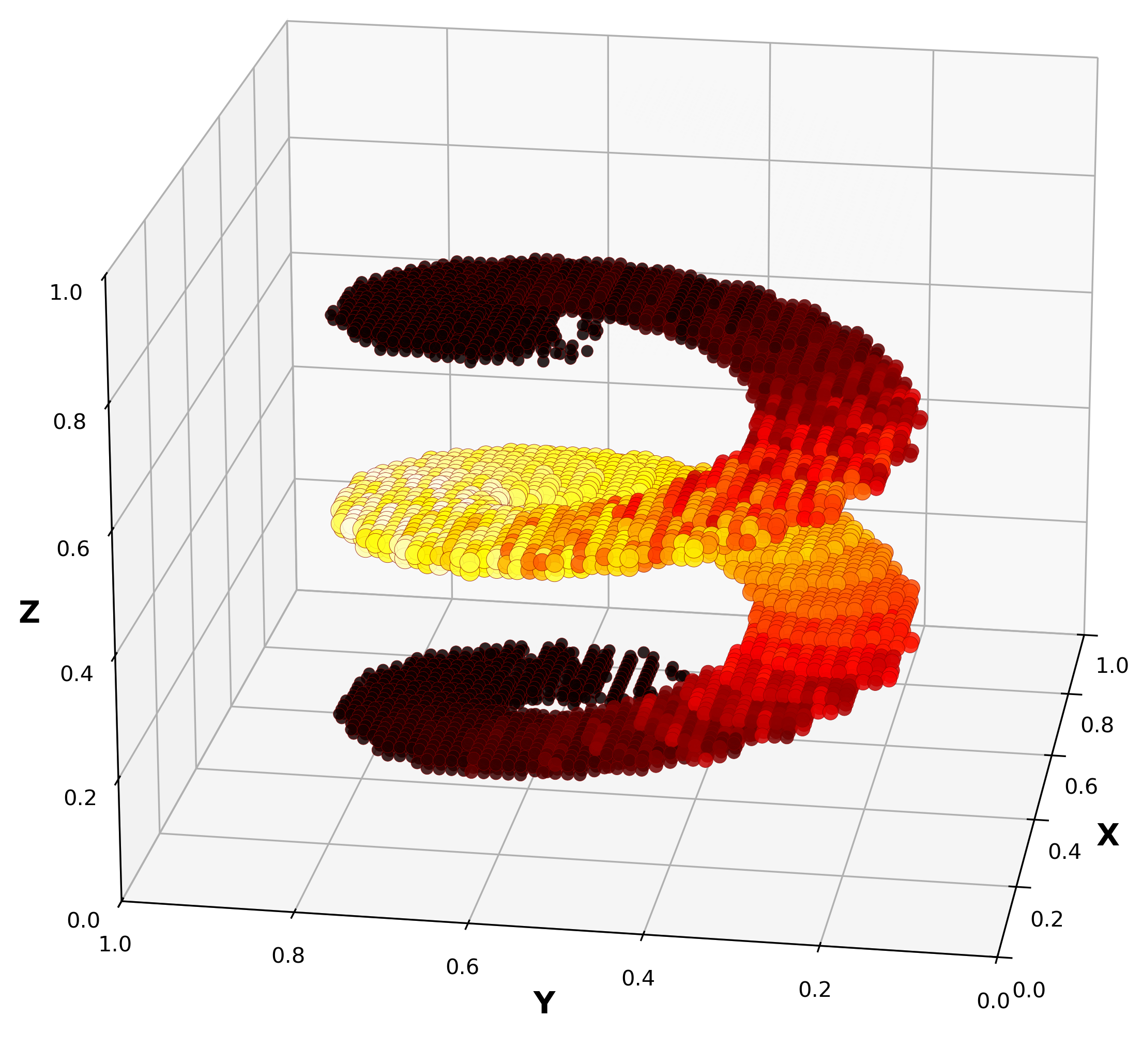}
        \caption{$t=0.4$}
    \end{subfigure}\hfill
    \begin{subfigure}{0.48\columnwidth}
        \centering
        \includegraphics[width=\linewidth]{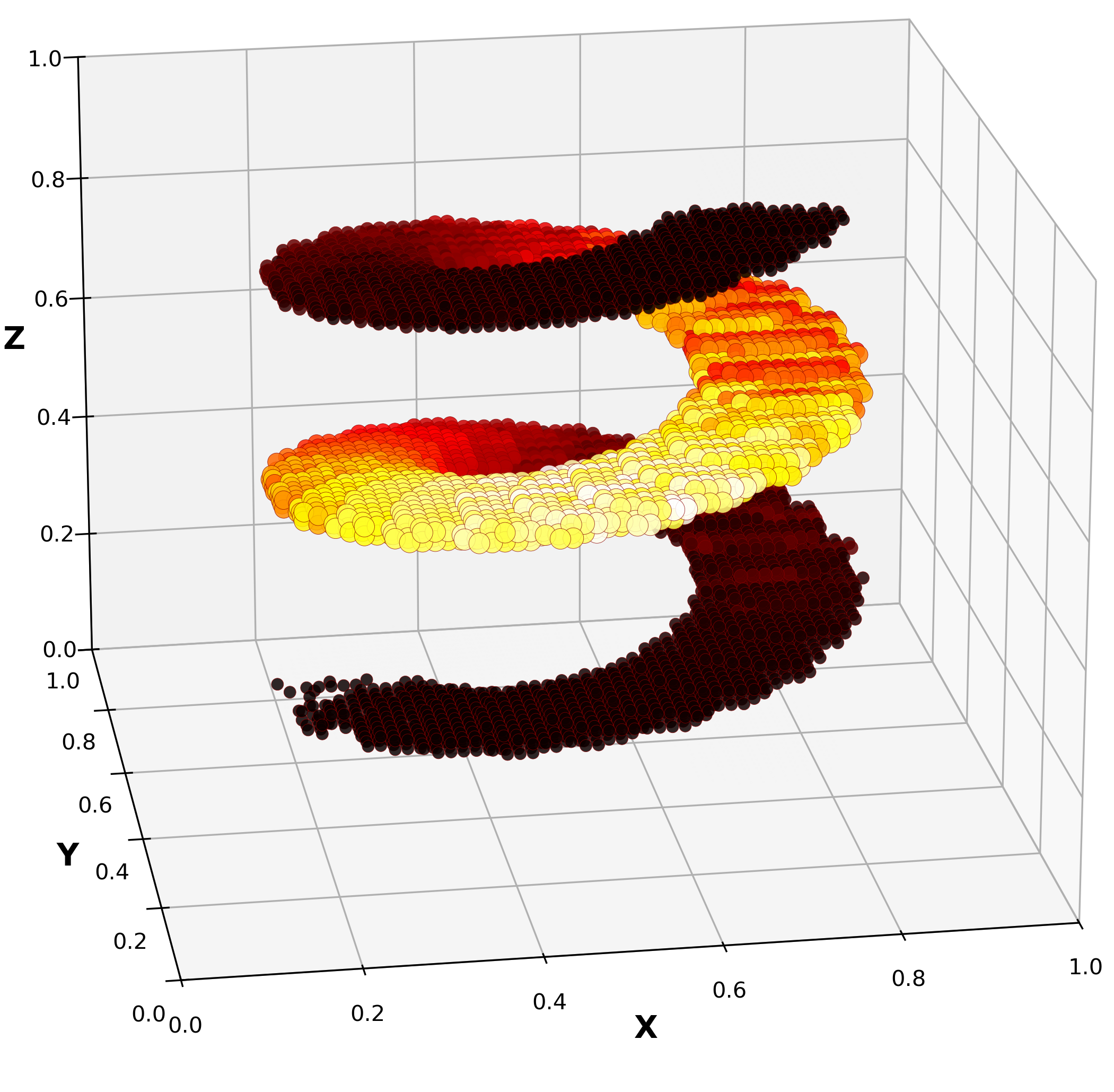}
        \caption{$t=0.6$}
    \end{subfigure}\\[4pt]
    \begin{subfigure}{0.48\columnwidth}
        \centering
        \includegraphics[width=\linewidth]{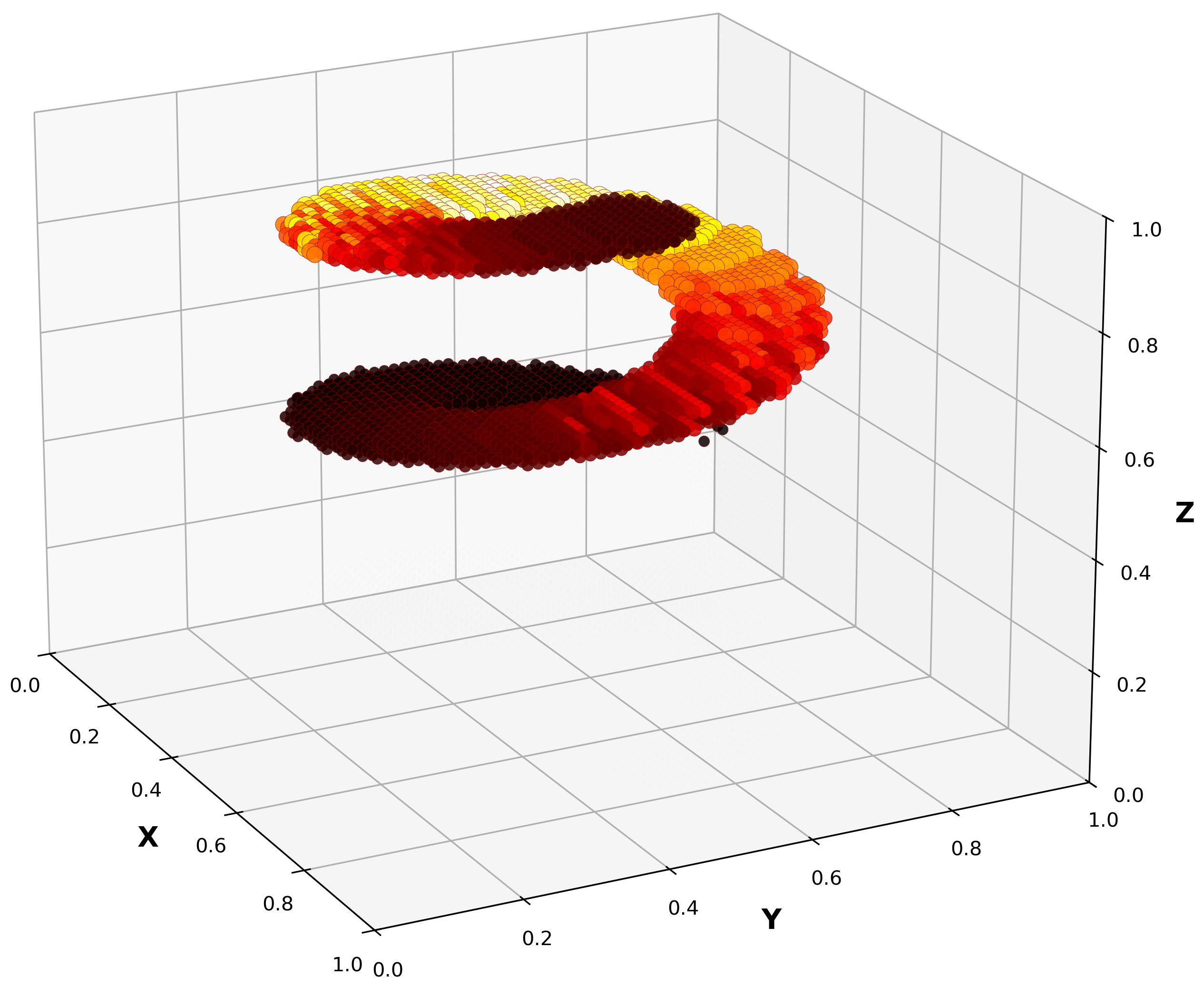}
        \caption{$t=0.8$}
    \end{subfigure}\hfill
    \begin{subfigure}{0.48\columnwidth}
        \centering
        \includegraphics[width=\linewidth]{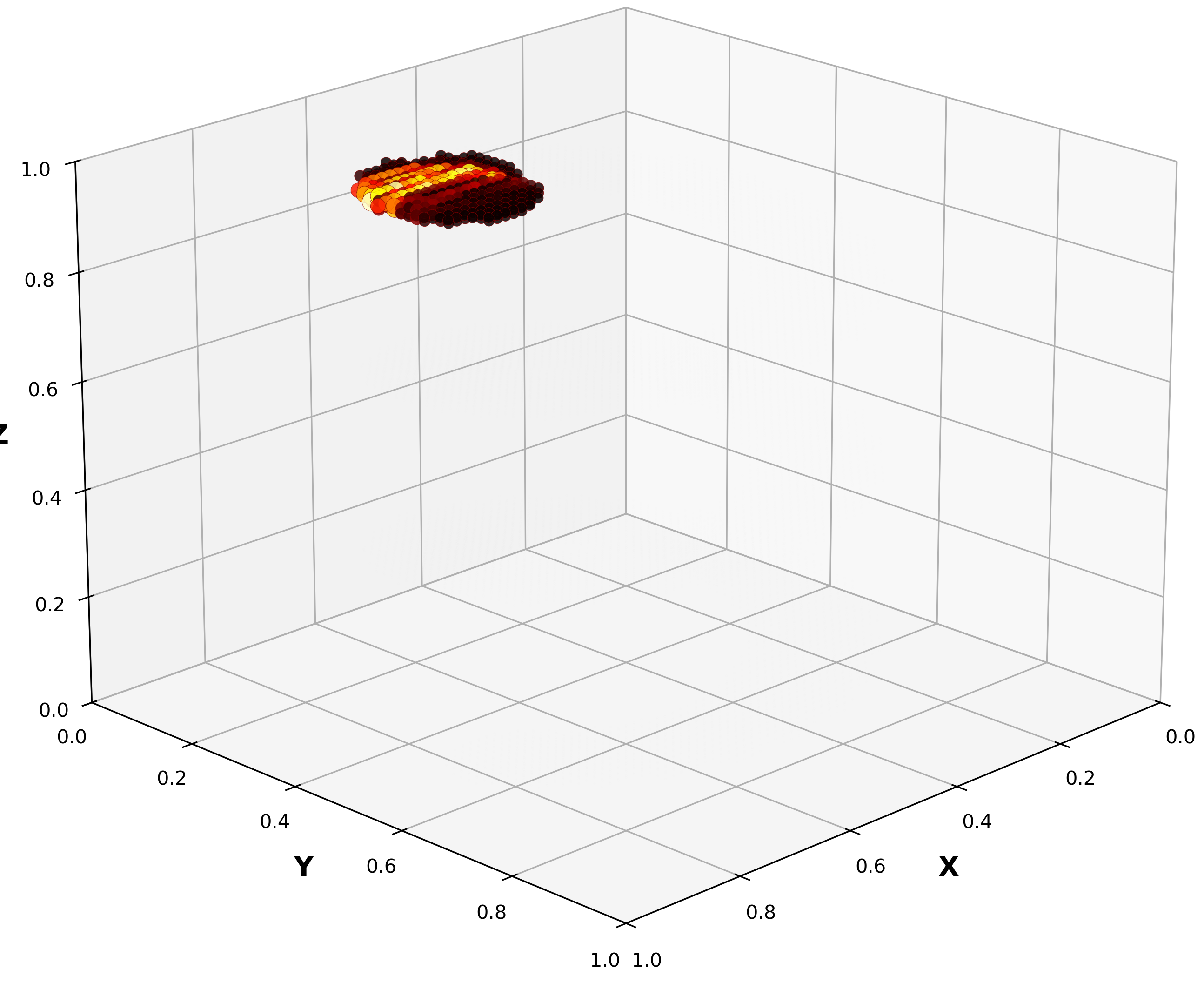}\caption{$t=1$}
    \end{subfigure}
    \caption{RSBP with prior drift: 3D temporal evolution of density $\rho(t,x)$ along the spiral channel.}
    \label{fig:drift_density_evolution_3d}
\end{figure}

Figure~\ref{fig:drift_particle_trajectories}(a) depicts the prior drift field. The drift exhibits a clear upward helical trend, pre-aligning particle motion along the spiral channel. Figure~\ref{fig:drift_density_evolution_3d} shows the complete 3D density evolution; the drift effectively guides the transport so that the optimal control $u^*$ focuses on fine-grained adjustments rather than driving the entire motion. The optimal control $u^*$ enters the reflecting SDE
\begin{equation}\label{sde_part}
dX_t = \bigl(v(X_t) + u^*(X_t,t)\bigr)\,dt + \sqrt{\varepsilon}\,dW_t + n(X_t)\,d\gamma_t
\end{equation}
as a feedback law. Figure~\ref{fig:drift_particle_trajectories}(b) shows the evolution of 10 particles under an Euler--Maruyama discretization of~\eqref{sde_part}. The controlled trajectories successfully navigate the spiral from entrance to exit.

We verify cost convergence of the scheme under spatial and temporal refinement. The spatial study fixes $K = 400$ and varies $N \in \{60, 70, 80, 90\}$; the temporal study fixes $N = 90$ and varies $K \in \{200, 300, 400\}$. Errors are computed relative to the finest discretization ($N=90$ or $K=400$).

\begin{table}[h]
\centering
\caption{3D spiral maze: spatial (left, $K=400$) and temporal (right, $N=90$) convergence of the cost $J$.}
\label{tab:convergence_energy}
\renewcommand{\arraystretch}{1.1}
\begin{tabular}{cccc}
\hline
$N$ & $J$ & $|\Delta J|$ & Rel.\ err.\ (\%) \\
\hline
60 & 3.887 & ---         & 0.95 \\
70 & 3.934 & $4.6\times10^{-2}$ & 0.23 \\
80 & 3.930 & $3.6\times10^{-3}$ & 0.14 \\
90 & 3.924 & $5.5\times10^{-3}$ & 0.00 \\
\hline
\end{tabular}
\hfill
\begin{tabular}{cccc}
\hline
$K$ & $J$ & $|\Delta J|$ & Rel.\ err.\ (\%) \\
\hline
200 & 3.531 & $2.2\times10^{-1}$ & 10.03 \\
300 & 3.784 & $1.1\times10^{-1}$ &  3.59 \\
400 & 3.924 & $6.2\times10^{-2}$ &  0.00 \\
\hline
\end{tabular}
\end{table}

\begin{table}[h]
\centering
\footnotesize
\caption{3D spiral maze: total mass $|\int_\Omega\rho\,dx - 1|$ at five time points. Spatial sweep at fixed $K=400$; temporal sweep at fixed $N=90$.}
\label{tab:convergence_mass}
\renewcommand{\arraystretch}{1.05}
\setlength{\tabcolsep}{3.5pt}
\begin{tabular}{cccc@{\quad}ccc}
\hline
& \multicolumn{3}{c}{Spatial sweep ($K=400$)} & \multicolumn{3}{c}{Temporal sweep ($N=90$)} \\
\cline{2-4}\cline{5-7}
$t$ & $N=70$ & $N=80$ & $N=90^*$ & $K=200$ & $K=300$ & $K=400$ \\
\hline
0.00 & $1.34\text{e}{-2}$ & $1.21\text{e}{-2}$ & $1.24\text{e}{-2}$ & $1.24\text{e}{-2}$ & $1.24\text{e}{-2}$ & $1.24\text{e}{-2}$ \\
0.25 & $1.56\text{e}{-2}$ & $1.85\text{e}{-2}$ & $1.78\text{e}{-2}$ & $1.85\text{e}{-2}$ & $1.80\text{e}{-2}$ & $1.78\text{e}{-2}$ \\
0.50 & $1.74\text{e}{-2}$ & $1.90\text{e}{-2}$ & $1.79\text{e}{-2}$ & $1.88\text{e}{-2}$ & $1.82\text{e}{-2}$ & $1.79\text{e}{-2}$ \\
0.75 & $1.50\text{e}{-2}$ & $1.49\text{e}{-2}$ & $1.63\text{e}{-2}$ & $1.72\text{e}{-2}$ & $1.66\text{e}{-2}$ & $1.63\text{e}{-2}$ \\
1.00 & $0$ & $0$ & $0$ & $0$ & $0$ & $0$ \\
\hline
\end{tabular}
\end{table}

Table~\ref{tab:convergence_energy} shows that $J$ varies by less than $1\%$ across all mesh resolutions, confirming mesh independence from $N = 70$ onwards, and converges monotonically with $K$. Table~\ref{tab:convergence_mass} shows mass deviations remaining below $2\%$ throughout, with the spatial and temporal sweeps converging towards the same reference values.


\bibliographystyle{IEEEtran}
\bibliography{bibliography_abbreviated}

\end{document}